\documentclass[leqno,10pt]{amsart}

\usepackage{amsmath,amssymb,amsthm}
\usepackage{epsfig}

\newtheorem{theorem}{Theorem}[section]

\newtheorem{lemma}[theorem]{Lemma}

\theoremstyle{remark}
  \newtheorem{remark}{Remark}[section]

\theoremstyle{definition}
  \newtheorem{definition}{Definition}[section]

\begin{document}

\title[Generalized splines in $\mathbb{R}^{n}$ and optimal
control]{Generalized splines in
$\mathbf{\mathbb{R}^{\lowercase{n}}}$\\
and optimal control}

\author{Rui~C.~Rodrigues, Delfim~F.~M.~Torres}

\thanks{This research was partially presented, as an oral
communication, at the Second Junior European Meeting on ``Control
Theory and Stabilization'', Dipartimento di Matematica del
Politecnico di Torino, Torino, Italy, 3-5 December 2003. To appear
on Rend. Sem. Mat. Univ. Pol. Torino, Vol. 64 (2006) No.1.}

\date{June 14, 2005}


\begin{abstract}
We give a new time-dependent definition of spline curves in
$\mathbb{R}^n$, which extends a recent definition of vector-valued
splines introduced by Rodrigues and Silva Leite for the
time-independent case. Previous results are based on a variational
approach, with lengthy arguments, which do not cover the
non-autonomous situation. We show that the previous results are a
consequence of the Pontryagin maximum principle, and are easily
generalized using the methods of optimal control. Main result
asserts that vector-valued splines are related to the Pontryagin
extremals of a non-autonomous linear-quadratic optimal control
problem.
\end{abstract}


\maketitle


\section{Introduction.}
Polynomial splines have been extensively used in several applied
areas of mathematics such as computer graphics and approximation
theory. Since the early 90's, they have been used in control
theory, associated to problems of aircraft control and path
planning of mechanical systems. These applications originated the
extension of classical spline functions to other contexts such as
Riemannian manifolds, Lie groups, etc.

Another line of research started with the definition of spline
functions which are not polynomial splines. One of the first
generalizations in this direction are the so called \emph{scalar
generalized splines}, which were introduced in the 50's by
Ahlberg, Nilson and Walsh \cite{ahlberg67}. The connection between
scalar generalized splines and optimal control was established
between 1995 and 1999. It turns out that splines are much more
than a tool to be used in control theory. They are intrinsic to
optimal control problems and appear naturally as minimizers of
certain problems \cite{clyde95,rui-leite-simoes99}.

Recently, this connection between minimality and splines was
extended to a new class of spline functions in arbitrary
dimensional Euclidean spaces \cite{rui-leite02-2}. This was
accomplished by variational arguments and a more general
time-invariant optimal control problem. Here, using tools from
optimal control, we go a step further. We consider a class of
classical linear-quadratic optimal control problems, which are not
necessarily time-invariant, and recover, as corollaries, the
previous results.


\section{Background.}
\label{sec:BackRD}

In this section we give an account of scalar generalized splines,
its connection to optimal control, and collect all the necessary
results to be used in the sequel.


\subsection{Scalar generalized splines.}
\label{subsec:sgs}

Generalized splines were first introduced in the late 50's by
Ahlberg, Nilson and Walsh \cite{ahlberg67}. Consider the linear
differential operator of order $p\in\mathbb{N}$
\[L = D^p \cdot + \,a_{p-1}(t)D^{p-1}\cdot + \cdots + a_1(t)D\cdot + a_0(t)\cdot\]
where each $a_k(t)$, $k = 0, 1, \ldots, p-1$, is a real
$C^{\,p}$-smooth function in $[a,b]$. The operator $L$ is acting
on the space $C^{\,m}[a,b]$ of real functions defined in $[a,b]$.
Its adjoint is defined by
\[L^* = (-1)^p D^p\cdot + \,(-1)^{p-1}D^{p-1}(a_{p-1}(t)\cdot)
+ \cdots - D(a_1(t)\cdot) + a_0(t)\cdot \,.\] $L^*$ is also acting
on $C^{\,m}[a,b]$ and the scalar product for which it is computed
is given by
\[\langle x_1 , x_2 \rangle = \int_{a}^{b} x_1(t)x_2(t) \,dt\,.\]
Let $\Delta : a = t_0 < t_1 < \ldots < t_m = b$, $m \in
\mathbb{N}$, be a partition of $[a,b]$, $\Omega$ be the family of
real $C^{\,2p-2}$-smooth functions in $[a,b]$ which are
$C^{\,2p}$-smooth in each interval $[t_{i},t_{i+1} ]$, $i = 0, 1,
\ldots, m-1$ and $f\in\Omega$.

\begin{definition}
\label{def:sg} The function $s : [a,b] \rightarrow \mathbb{R}$ is
an interpolating generalized spline of $f$ associated to $\Delta$
and $L$, if $s \in \Omega$, $s$ is a solution of the differential
equation $L^*Lx=0$ in each interval $[t_{i},t_{i+1}]$, $i = 0, 1,
\ldots, m-1$, and $s(t)=f(t)$ on $\Delta$.
\end{definition}

\begin{definition}
\label{def:sg-tipoI} An interpolating generalized spline of $f$ is
of type I if it is such that $s^{(k)}(t_0) = f^{(k)}(t_0)$ and
$s^{(k)}(t_m)=f^{(k)}(t_m)$, for $k = 1, 2, \ldots ,p-1$.
\end{definition}

The function $f$ is usually omitted from the previous definitions.
Instead, one has to demand that, in Definition~\ref{def:sg},
function $s$ fulfills the interpolation condition $s(t_i) = s_i$,
where $s_i$, $i = 0,1,\ldots,m$, are given real numbers and, in
Definition~\ref{def:sg-tipoI}, that $s$ fulfills the boundary
conditions $s^{(k)}(t_0) = \eta_0^k$ and $s^{(k)}(t_m) = \eta_m^k$
where $\eta_0^k$, $\eta_m^k$, $k = 1, 2, \ldots, p-1$, are
prescribed real numbers. Then, we just say that $s$ is a
\emph{generalized spline} of type $I$. The next statement collects
several results about generalized splines of type $I$ which can be
found in \cite{ahlberg67}.
\begin{theorem}[\cite{ahlberg67}]
There exists, for each set of boundary and interpolation
conditions, a unique generalized spline of type $I$ associated
with the differential operator $L$ and the partition $\Delta$.
Moreover, this generalized spline is the unique minimizer of the
functional
\[\int_{a}^{b} (Lg)^2 \,dt\]
among all the functions $g \in \Omega $ that fulfill the same
boundary and interpolation conditions.
\end{theorem}

\begin{remark}
There are other types of boundary conditions, described in the
literature, that also ensure the existence and uniqueness of the
corresponding generalized spline.
\end{remark}

We now give two examples for constant coefficient operators: an
example of a cubic spline, and an example of a trigonometric
spline. Let $\Delta : 0 < 1/4 < 1$ be the partition of the time
interval $[0,1]$; $s(t_0) = 3$, $s(t_1) = 1$ and $s(t_2) = 0$ be
the interpolation conditions; and $\dot{s}(t_0) = -1$,
$\dot{s}(t_2)=1$ be the boundary conditions. We first consider the
operator $L=D^2$. The resulting spline of type I is a
$C^{\,2}$-smooth function in $[0,1]$ such that
$s(t)=c_{1i}+c_{2i}t+c_{3i}t^2+c_{4i}t^3$ in each
$[t_{i},t_{i+1}]$ where $c_{1i}$, $c_{2i}$, $c_{3i}$, $c_{4i}$ are
real constants to be found. This is the classical cubic spline.
Considering $L=D^2+144$, the resulting spline of type I is also
$C^{\,2}$-smooth in $[0,1]$ so that
$s(t)=(c_{1i}+c_{2i}t)\cos{(12t)}+(c_{3i}+c_{4i}t)\sin{(12t)}$ in
each $[t_{i},t_{i+1}]$.

\medskip

\noindent
\begin{minipage}[b]{.46\linewidth}
\begin{center} \epsfig{figure = 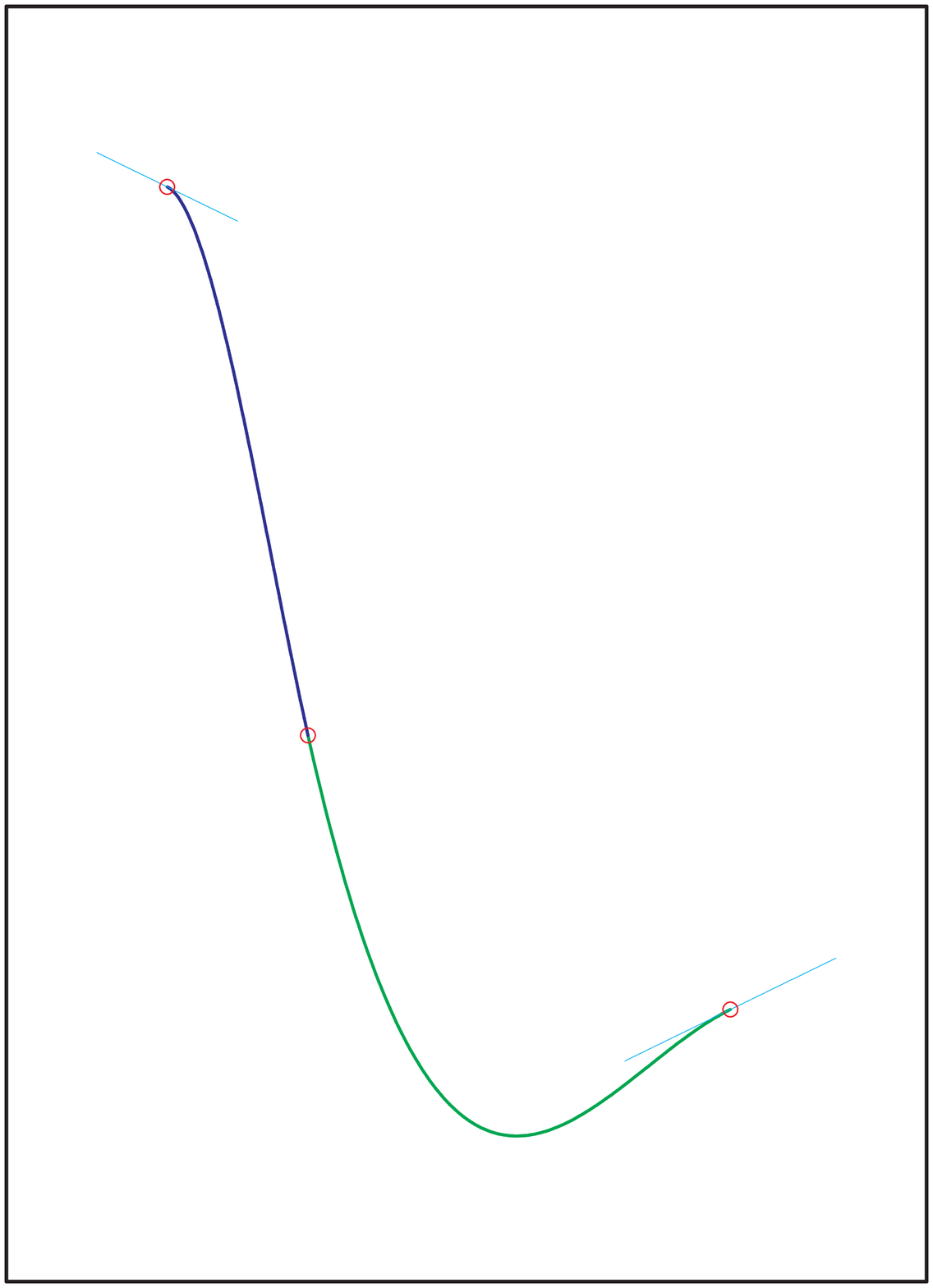, height = 5cm, width = 5.5cm} \\
\mbox{} \\ Cubic spline \end{center}
\end{minipage} \hfill
\begin{minipage}[b]{.46\linewidth}
\begin{center} \epsfig{figure = 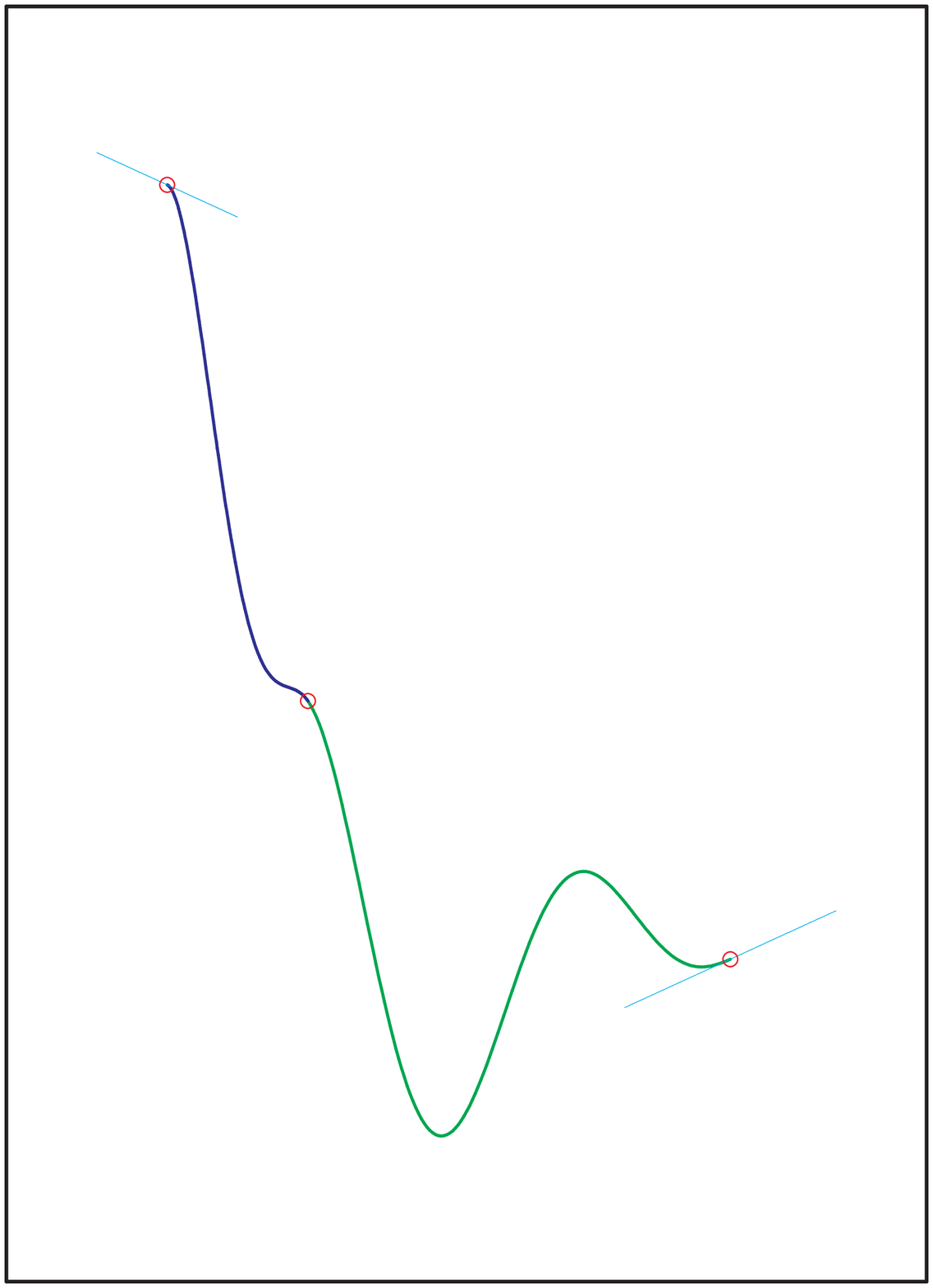, height = 5cm, width = 5.5cm} \\
\mbox{} \\  Trigonometric spline \end{center}
\end{minipage}

\medskip

The most immediate generalization of scalar splines to curves in
$\mathbb{R}^n$ is achieved by simply considering vector functions
$g : [a,b] \rightarrow \mathbb{R}^n$, the same operator $L$ as
before, and adapted interpolation conditions, boundary conditions,
and set $\Omega$. It is obvious that each component of the
resulting spline will be a scalar generalized spline, and
therefore such a spline curve will always minimize the functional
\[\int_{a}^{b} \langle Lg, Lg \rangle \,dt \, ,\]
where $\langle \cdot \, , \cdot \rangle$ stands for the Euclidean
inner product, among all functions in $\Omega$ that fulfill the
same boundary and interpolation conditions. As we shall see, from
an optimal control perspective such a trivial generalization is
not the natural way of extending scalar-splines to vector-valued
splines.


\subsection{Scalar generalized splines and optimal control.}
Since the early nineties, in order to deal with applied problems
from Robotics, there has been an increasing interest to combine
spline curves and integral cost problems associated with linear
control systems. Among theoretical developments, it was found that
scalar generalized splines are minimizers of a simple optimal
control problem with a linear time-invariant control system and a
single control (see~\cite{clyde95} and \cite{rui-leite-simoes99}).
This discovery is of crucial importance, because it introduces a
new perspective to the subject: scalar spline functions are better
viewed as a consequence of the search for an optimal control,
rather than a postulate imposed a priori in order to solve
particular classes of problems. Given its importance, we summarize
the main result here. Consider the following autonomous
linear-quadratic optimal control problem:
\begin{equation}
\label{eq:SGSOC}
\begin{split}
& \min_{u} \quad \int_{a}^{b} u^2 \,dt \\
\text{subject to} \\
& \dot{x} = Ax + Bu \\
& x(t_0)=x_0, \quad x(t_{m})=x_{m} \\
& x^1(t_i) = \alpha_i, \qquad i = 1,2, \ldots, m-1 \, ,
\end{split}
\end{equation}
where $a = t_0 < t_1 < \cdots < t_{m-1} < t_{m} = b$, $x^1$ is the
first component of the state vector, $\alpha_i\in\mathbb{R}$, $u$
is a scalar function which is $C^{\,n-2}$-smooth in $[a,b]$ and
$C^{\,n}$-smooth in each interval $[t_i,t_{i+1}]$. Let us assume
that the state space is $\mathbb{R}^n$ and that the state vector
is a $C^{\,2n-2}$-smooth function in $[a,b]$ which is also
$C^{\,2n}$-smooth in each interval $[t_i,t_{i+1}]$.

\begin{theorem}
If the control system $\dot{x} = Ax + Bu$ of problem
\eqref{eq:SGSOC} is completely state controllable with matrices
$A$ and $B$ in the canonical form {\small
\[A=\left(
\begin{array}{cccc}
0 & 1 & \cdots & 0 \\
\vdots & \ddots & \ddots & \vdots \\
0 & \cdots & 0 & 1 \\
a_0 & a_1 & \cdots & a_{n-1}
\end{array}\right), \qquad
B=\left(\begin{array}{c} 0 \\ \vdots \\ 0 \\ 1
\end{array}\right)\]}
for given real numbers $a_j$, $j=0,1,\ldots,n-1$, then the optimal
control problem \eqref{eq:SGSOC} has always a unique solution with
the first component of the optimal state vector being a
generalized spline of type $\rm I$ associated to the constant
coefficient differential operator $L = D^n - \,a_{n-1}D^{n-1} -
\cdots - a_1D - a_0$.
\end{theorem}

If the first component of the optimal state vector of problem
\eqref{eq:SGSOC} is a scalar generalized spline, the following
questions come immediately to our mind: What can be said about the
minimizing state trajectory of the optimal control problem? Is it
some sort of a generalized spline in $\mathbb{R}^n$? The answer to
these questions leads us (see Definition~\ref{def:VGS:TD}) to a
new time-dependent definition of ge\-ne\-ra\-lized spline in
$\mathbb{R}^n$, which is the main contribution of the present
paper.


\subsection{Pontryagin's maximum principle, existence, and regularity.}
The general problem of optimal control can be defined, in Lagrange
form, as follows:
\begin{gather}
\min_{\left(x(\cdot),u(\cdot)\right)} \quad I[x(\cdot),u(\cdot)] =
\int_{t_a}^{t_b}
\mathcal{L}\left(t,x(t),u(t)\right)\,\mathrm{d}t \notag \\
\dot{x}\left(t\right) =
\varphi\left(t,x\left(t\right),u\left(t\right)\right) \label{eq:prb:Pgen}\\
\left(x(t_a),x(t_b)\right) = (\alpha,\beta)  \notag \\
x\left(\cdot\right) \in W_{1,1}\left(\left[t_a,t_b\right];
\,\mathbb{R}^{n}\right) \, \notag \\
u\left( \cdot \right) \in \mathcal{U}\left(\left[t_a,t_b\right];\,
\Omega \subseteq \mathbb{R}^{r}\right) \, . \notag
\end{gather}
We assume that $\mathcal{L}:\left[t_a,t_b\right] \times
\mathbb{R}^{n}\times \mathbb{R}^r\rightarrow \mathbb{R}$ and
$\varphi: \left[t_a,t_b\right] \times \mathbb{R}^{n} \times
\mathbb{R}^r \rightarrow \mathbb{R}^{n}$ are $C^1$-smooth
functions with respect to all arguments, and that the boundary
conditions, together with the class of control functions
$\mathcal{U}$, are given. The standard method to solve
\eqref{eq:prb:Pgen} is usually based on the deductive approach:
(i) a solution exists for the problem; (ii) the necessary
conditions are applicable, and they identify certain candidates
(so called \emph{extremals}); (iii) subsequent elimination (if
necessary) identifies the solution (or solutions). We are
interested in the case where there are no restrictions on the
control variables: $\Omega = \mathbb{R}^{r}$. The unrestricted
case poses many difficulties, and the problem turns out to be a
difficult one, even in special situations. As we explain next,
most part of difficulties appear in the application of steps (i)
and (ii).

The first general answer to (i) was given by A. F. Filippov in
1959 \cite{MR22:13373}, assuming the admissible controls to be
integrable ($\mathcal{U} = L_1$), and the control set $\Omega$ to
be compact. As far as we assume $\Omega$ to be a noncompact set,
Filippov's theorem does not apply. To solve the existence problem,
we make use of the following theorem (see \cite{MR85c:49001}).

\begin{theorem}[``Tonelli'' Existence Theorem for \eqref{eq:prb:Pgen}]
\label{Th:Ex:Tonelli} Problem \eqref{eq:prb:Pgen} has a minimizer
$\left(\tilde{x}(\cdot),\tilde{u}(\cdot)\right)$ with
$\tilde{u}(\cdot) \in L_1\left([t_a,t_b];\,\mathbb{R}^r\right)$,
provided there exists at least one admissible pair, and the
following convexity and coercivity conditions hold:
\begin{itemize}
\item \emph{(convexity)} Functions $\mathcal{L}(t,x,\cdot)$ and
$\varphi(t,x,\cdot)$  are convex for all $\left(t,x\right)$;
\item \emph{(coercivity)} There exists a function
$\theta : \mathbb{R}_{0}^{+} \rightarrow \mathbb{R}$, bounded
below, such that
\begin{gather*}
\mathcal{L}(t,x,u) \ge
\theta\left(\left\|\varphi(t,x,u)\right\|\right)
\quad \text{for all } (t,x,u) \text{;} \label{e:coercivity1}\\
\lim_{r \rightarrow + \infty} \frac{\theta(r)}{r} = +
\infty\text{;} \label{e:coercivity2} \\
\lim_{\left\|u\right\| \rightarrow + \infty}
\left\|\varphi(t,x,u)\right\| = + \infty \quad \text{for all }
(t,x)\, . \label{e:coercivity3}
\end{gather*}
\end{itemize}
\end{theorem}
\begin{remark}
For the definition of convexity of $\mathcal{L}(t,x,\cdot)$ and
$\varphi(t,x,\cdot)$ see \cite{MR85c:49001}. In the case $\varphi
= u$ one has the fundamental problem of the calculus of
variations, and we get from Theorem~\ref{Th:Ex:Tonelli} the
classical Tonelli existence theorem.
\end{remark}

Step (ii) is addressed by the Pontryagin Maximum Principle
\cite{pontryagin62}.
\begin{theorem}[Pontryagin Maximum Principle]
\label{Teo:PMP} If $\left(x(\cdot),\,u(\cdot)\right)$ is a
minimizer of \eqref{eq:prb:Pgen} and $u(\cdot)$ is essentially
bounded, $u(\cdot) \in L_{\infty}$, then there exists
$\left(\psi_{0},\,\psi(\cdot)\right) \ne 0$, $\psi_0 \le 0$,
$\psi(\cdot) \in W_{1,1}^n$, such that the quadruple
$\left(x(\cdot),u(\cdot),\psi_{0},\psi(\cdot)\right)$ is a
Pontryagin extremal: it satisfies
\begin{itemize}
\item the \emph{Hamiltonian system}
\begin{equation}
\label{PMP:HamSyst}
\begin{cases}
\dot{x}=\dfrac{\partial H}{\partial\psi} \, ,\\
\dot{\psi}=-\dfrac{\partial H}{\partial x} \, ;
\end{cases}
\end{equation}
\item the \emph{maximality condition}
\begin{equation}
\label{PMP:maxCond} H\left(t,x(t),u(t),\psi_{0},\psi(t)\right)
=\max_{v\in\mathbb{R}^r}H\left(t,x(t),v,\psi_{0},\psi(t)\right)\,
;
\end{equation}
\end{itemize}
with the Hamiltonian
\begin{equation}
\label{PMP:Ham}
H(t,x,u,\psi_0,\psi)=\psi_{0}\,\mathcal{L}\left(t,x,u\right) +
\langle \psi,\varphi\left(t,x,u\right) \rangle \, .
\end{equation}
\end{theorem}

\begin{definition}
A Pontryagin extremal
$\left(x(\cdot),u(\cdot),\psi_{0},\psi(\cdot)\right)$ is said to
be abnormal when $\psi_0$ is equal to zero, and normal otherwise.
\end{definition}

The existence is assured in the class of integrable controls
($\mathcal{U} = L_{1}$), while the formulation of the Pontryagin
maximum principle assume the optimal controls to be essentially
bounded ($\mathcal{U} = L_{\infty} \subset L_{1}$). For minimizers
predicted by existence theory, Theorem~\ref{Teo:PMP} may fail to
be valid, because the values of optimal controls can be unbounded.
This is a possibility even for very simple instances of problem
\eqref{eq:prb:Pgen}: \textrm{e.g.} $\mathcal{L}$ a polynomial and
$\varphi$ linear. One such example can be found in
\cite{MR86k:49002}: the problem
\begin{gather}
\min \quad \int_{0}^{1}\left( \left( x^{3}-t^{2}\right) ^{2}\,
u^{14}
+\varepsilon \,u^{2}\right) \,dt \notag \\
\dot{x}\left(t\right) =u\left(t\right) \label{eq:prb:BM}\\
x\left(0\right) =0\, , \quad x\left(1\right) =k \, , \notag
\end{gather}
satisfies all the hypotheses of Theorem~\ref{Th:Ex:Tonelli}; it
can be proved (see \cite{MR85m:49051}) that for certain choices of
constants $k$ and $\varepsilon$ there exists a unique optimal
control $u\left( t\right) =k\,t^{-1/3}$; but Theorem~\ref{Teo:PMP}
(Pontryagin maximum principle) is not satisfied since
$\dot{\psi}(t) =
\mathcal{L}_{x}\left(t,x\left(t\right),\dot{x}\left(t\right)\right)
= c\,t^{-4/3}$ is not integrable ($\psi(\cdot)$ is not an
absolutely continuous function).

In order to apply the deductive method (i)--(iii) one needs to
close the gap between the hypotheses of existence and necessary
optimality conditions. For that, conditions beyond those of
convexity and coercivity, assuring solutions $\tilde{u}(\cdot)$ to
be in $L_{\infty}$ and not only in $L_{1}$, must apply. To exclude
the possibility of bad behavior that occurs for \eqref{eq:prb:BM},
we will focus our attention to problem \eqref{eq:prb:Pgen} with
Lagrangian $\mathcal{L}$ and function $\varphi$ given by
\begin{equation}
\begin{split}
\label{eq:bigNoGap} \mathcal{L}(t,x,u) = \langle B(t)u,B(t)u \rangle \, ,
\\
\varphi(t,x,u) = A(t)x + B(t)u \, ,
\end{split}
\end{equation}
under the hypothesis
\begin{description}
\item[(H1)] $B(t)$ is a square matrix with full rank for all $t$;
\item[(H2)] the dynamical control system $\dot{x}(t) = A(t)x(t) +
B(t)u(t)$ is completely state controllable; \item[(H3)] $t
\rightarrow A(t)$ and $t \rightarrow B(t)$ are $C^{1}$-smooth
functions.
\end{description}
Roughly speaking, this gives the biggest class of optimal control
problems which generalize \eqref{eq:SGSOC} in a natural way; do
not admit abnormal extremals (see the next remark); and for which
the gap between existence and necessary optimality conditions is
automatically closed.

\begin{remark}
Since there is no constraint on the control, singular
trajectories are exactly projections of abnormal extremals. But
due to the assumption on the linear system (it is supposed to be
completely state controllable), there is no singular trajectory,
and thus the optimal control problem has no abnormal extremals.
\end{remark}

\begin{theorem}[Boundedness of optimal controls \cite{delfim03}]
\label{teo:regDT} Under the hypotheses of the existence
Theorem~\ref{Th:Ex:Tonelli}, if there exist constants $c>0$ and
$k$ such that
\begin{equation*}
\begin{split}
\left|\frac{\partial \mathcal{L}}{\partial t}\right| &\le c
\left|\mathcal{L}\right| + k \, , \quad \left\|\frac{\partial
\mathcal{L}}{\partial x}\right\| \le
c \left|\mathcal{L}\right| + k \, ,\\
\left\|\frac{\partial \varphi}{\partial t}\right\| &\le c
\left\|\varphi\right\| + k \, , \quad \left\|\frac{\partial
\varphi_i}{\partial x}\right\| \le c \left|\varphi_i\right| + k
\quad (i = 1,\,\ldots,\,n) \, \text{;}
\end{split}
\end{equation*}
then all minimizers of \eqref{eq:prb:Pgen} satisfy the Pontryagin
maximum principle.
\end{theorem}

It is a simple exercise to see that with $\mathcal{L}$ and
$\varphi$ defined by \eqref{eq:bigNoGap}, hypotheses (H1) and (H3)
imply all the conditions of Theorems~\ref{Th:Ex:Tonelli} and
\ref{teo:regDT}.


\section{Main results.}
We are interested in the following non-autonomous linear-quadratic
optimal control problem:
\begin{equation}
\tag{$P$}
\begin{split}
& \min_{u(\cdot)} \quad J\left[u(\cdot)\right]
= \int_{a}^{b} \langle B(t)u(t),B(t)u(t) \rangle \,dt \\
\text{subject to} \\
& \dot{x}(t) = A(t)x(t) + B(t)u(t) \\
& x(t_i)=x_i, \qquad i = 0,1,\ldots,m \, ,
\end{split}
\end{equation}
for a given partition $a = t_0 < t_1 < \cdots < t_{m-1} < t_{m} = b$
and fixed $x_i \in \mathbb{R}^n$. The control $u : [a,b] \rightarrow
\mathbb{R}^{n}$ is unrestricted; the state function $x : [a, b]
\rightarrow \mathbb{R}^n$ is an absolutely continuous function;
$A(t)$ and $B(t)$ are $n \times n$ matrices and $B(t)$ is
nonsingular. We find the minimizer of $(P)$ by solving $(P_i)$, $i =
0,1,\ldots,m-1$, in the interval $[t_{i},t_{i+1}]$:
\begin{equation}
\tag{$P_i$}
\begin{split}
& \min_{u(\cdot)} \quad J_i\left[u(\cdot)\right] = \int_{t_{i}}^{t_{i+1}} \langle B(t)u(t),B(t)u(t) \rangle \,dt \\
& \dot{x}(t) = A(t)x(t) + B(t)u(t) \\
& x(t_{i})=x_{i}, \qquad x(t_{i+1})=x_{i+1} \, .
\end{split}
\end{equation}
In order to guarantee the applicability of the Pontryagin maximum
principle, and the existence of a normal solution, hypotheses
(H1), (H2) and (H3) of previous section are in force. Under these
assumptions we can choose, without any loss of generality, $\psi_0
= -\frac{1}{2}$ in Theorem~\ref{Teo:PMP}. The Hamiltonian
\eqref{PMP:Ham} is then given by
\begin{equation*}
H(t,x,u,\psi)=-{\textstyle\frac{1}{2}}\,u^{\prime}B(t)^{\prime}B(t)u
+\psi^{\prime}(A(t)x+B(t)u) \, ,
\end{equation*}
where we use the symbol prime $^\prime$ to denote the transpose of
a given vector or matrix. The Hamiltonian system
\eqref{PMP:HamSyst} reduces to
\begin{equation}
\label{sistemahamiltoniano}
\begin{cases}
\dot{x}(t) = A(t)x(t)+B(t)u(t) \, ,\\
\dot{\psi}(t) = -A(t)^{\prime}\psi(t) \, ,
\end{cases}
\end{equation}
while from the maximality condition \eqref{PMP:maxCond} one
obtains
\begin{equation*}
B(t)^{\prime}(B(t)u(t) - \psi(t)) = 0 \, .
\end{equation*}
This equation implies that $\psi(t) = B(t)u(t)$ and hence
\begin{equation}
\label{eq:qs:ft} u(t) = B(t)^{-1}\psi(t)
\end{equation}
is the unique Pontryagin extremal control. Thus, due to
Theorem~\ref{teo:regDT}, $u$ given by \eqref{eq:qs:ft} must be
optimal.

From equation $\psi(t) = B(t)u(t)$ and from equation $\dot{\psi}(t)
= -A(t)^{\prime}\psi(t)$ of system \eqref{sistemahamiltoniano} we
get the matrix differential equation
\begin{equation}
\label{eq:mde} \frac{d}{dt}(B(t)u(t)) + A(t)^{\prime}B(t)u(t) = 0\,.
\end{equation}
Introducing the matrix differential operator $L=D-A(t)$, the control
system $\dot{x}(t) = A(t)x(t)+B(t)u(t)$ can be written as
\begin{equation}
\label{eq:mdeTrf1} Lx(t) = B(t)u(t)
\end{equation}
and equation \eqref{eq:mde} as
\begin{equation}
\label{eq:mdeTrf2} L^*B(t)u(t) = 0 \, ,
\end{equation}
where $L^* = -D-A(t)^{\prime}$ is the adjoint operator of $L$. From
\eqref{eq:mdeTrf1} and \eqref{eq:mdeTrf2} we conclude that the
minimizing state trajectory is a solution of the differential
equation
\begin{equation*}
L^*Lx(t) = 0 \, ,
\end{equation*}
which can be written as
\[\ddot{x}(t) + (A(t)^{\prime} - A(t))\,\dot{x}(t) - (A(t)^{\prime}A(t)+\dot{A}(t))\,x(t) = 0 \, .\]
We have just proved Lemma~\ref{lemaRCDT}.
\begin{lemma}
\label{lemaRCDT} Under hypotheses (H1)-(H3) the optimal control $u$
is, in each interval $[t_{i},t_{i+1}]$, $i=0,1,\ldots,m-1$, a
solution of the matrix differential equation $L^*B(t)u = 0$ with
$L^*=-D-A(t)^{\prime}$ the adjoint operator associated to the
operator $L=D-A(t)$. The corresponding optimal state trajectory $x$
is such that $L^*Lx=0$ in each interval $[t_{i},t_{i+1}]$.
\end{lemma}

An explicit expression for the optimal state trajectory and for
the optimal control can be obtained in terms of the \emph{state
transition matrix}. These results are stated in the following
Theorem. We refer the reader to \cite{brockett70} for the
definition, and properties, of the state transition matrix.

\begin{theorem}
\label{teo:Pre:2.12} The optimal state trajectory of problem $(P)$
has, in each interval $[t_{i},t_{i+1}]$, $i=0,1,\ldots,m-1$, the
explicit expression
\begin{equation}
\label{eq:ExpExpGRSL}
x(t)=\Phi(t,t_i)x_i+\left(\int_{t_i}^{t}\Phi(t,s)\Phi(t_i,s)^{\prime}
\,ds\right)S^{-1}(\Phi(t_i,t_{i+1})x_{i+1}-x_i)
\end{equation}
where $\Phi$ is the state transition matrix associated to
$\dot{x}=A(t)x$, and $S$ is the symmetric matrix given by
\[\int_{t_i}^{t_{i+1}} \! \Phi(t_i,s)\Phi(t_i,s)^{\prime}\,ds \, .\]
Furthermore, the optimal control of problem $(P)$ has, in each
interval $[t_{i},t_{i+1}]$, $i=0,1,\ldots,m-1$, the explicit
expression
\begin{equation}
\label{eq:optimalcontrol}
u(t)=B(t)^{-1}\,\Phi(t_i,t)^{\prime}\,S^{-1}
(\Phi(t_i,t_{i+1})\,x_{i+1}-x_i) \, .
\end{equation}
\end{theorem}

\begin{proof}
(Theorem~\ref{teo:Pre:2.12}) Since $\psi = B(t)u$, the Hamiltonian
system takes the form
\begin{equation}\begin{cases}
\dot{x}=A(t)x + \psi \, ,\\
\dot{\psi}=-A(t)^{\prime}\psi \, .
\end{cases}\label{sistemahamiltoniano2}\end{equation}
From equation $\dot{\psi}=-A(t)^{\prime}\psi$ we get
$\psi(t)=\Phi(t_i,t)^{\prime}\psi(t_i)$. The substitution of
$\psi$ in equation $\dot{x}=A(t)x + \psi$ of system
(\ref{sistemahamiltoniano2}) generates
$\dot{x}=A(t)x+\Phi(t_i,t)^{\prime}\psi(t_i)$. The solution of
this complete differential equation, with initial condition
$x(t_i)=x_i$, is given by
\begin{equation}
x(t) = \Phi(t,t_i)x_i +
\int_{t_i}^{t}\Phi(t,s)\Phi(t_i,s)^{\prime}\psi(t_i)\,ds.
\label{state}
\end{equation}
Now, we just have to find $\psi(t_i)$. Using the other initial
condition $x(t_{i+1})=x_{i+1}$ we get
\[\begin{split}
x_{i+1} &= \Phi(t_{i+1},t_i)x_i + \left(\int_{t_i}^{t_{i+1}}
\Phi(t_{i+1},s)\Phi(t_i,s)^{\prime}\,ds\right)\psi(t_i) \\
&= \Phi(t_{i+1},t_i)x_i + \Phi(t_{i+1},t_i)
\left(\int_{t_i}^{t_{i+1}}
\Phi(t_i,s)\Phi(t_i,s)^{\prime}\,ds\right)\psi(t_i) \, .
\end{split}\]
If we denote the symmetric matrix
\[\int_{t_i}^{t_{i+1}} \! \Phi(t_i,s)\Phi(t_i,s)^{\prime}\,ds\]
by $S(t_i,t_{i+1})$, or simply by $S$, we can write
\[\Phi(t_{i+1},t_i)^{-1}x_{i+1}-x_i=S\,\psi(t_i)
\Leftrightarrow \Phi(t_i,t_{i+1})x_{i+1}-x_i=S\,\psi(t_i).\] Since
matrix $S$ is always non-singular, we get
$\psi(t_i)=S^{-1}(\Phi(t_i,t_{i+1})x_{i+1} - x_i)$. Finally, from
(\ref{state}), we obtain the equality \eqref{eq:ExpExpGRSL}:
\[x(t)_{|t\in[t_i,t_{i+1}]} = \Phi(t,t_i)x_i + \left(\int_{t_i}^{t}
\Phi(t,s)\Phi(t_i,s)^{\prime}\,ds\right)
S^{-1}(\Phi(t_i,t_{i+1})x_{i+1} - x_i).\]%
The second part of the theorem is a direct consequence of equation
\eqref{eq:qs:ft}. From previous calculations we have
\[\psi(t) = \Phi(t_i,t)^{\prime}\,S^{-1}(\Phi(t_i,t_{i+1})x_{i+1} - x_i)\]
and thus, equality \eqref{eq:optimalcontrol} follows immediately.
\end{proof}

\begin{remark}
From the proof of Theorem~\ref{teo:Pre:2.12} it follows,
by direct calculations, that the optimal value for the
integral functional $J_i[\cdot]$ of problem $(P_i)$
is given by
\begin{equation*}
(\Phi(t_i,t_{i+1})x_{i+1} - x_i)^{\prime}\,S^{-1}\,
(\Phi(t_i,t_{i+1})x_{i+1} - x_i) \, .
\end{equation*}
\end{remark}
\begin{remark}
We have seen that in each interval $[t_{i},t_{i+1}]$,
$i=0,1,\ldots,m-1$, the optimal state trajectory of problem $(P)$ is
a solution of the matrix differential equation
\[\ddot{x}(t) + (A(t)^{\prime} - A(t))\,\dot{x}(t) - (A(t)^{\prime}A(t)+\dot{A}(t))\,x(t) = 0\]
which does not depend on the matrix $B(t)$. This is natural since we
can make  the substitution $u \mapsto v = B(t)u$ in the problem
$(P)$ and thus eliminate the presence of matrix $B(t)$ in all
further calculations.
\end{remark}
\begin{remark}
When problem $(P)$ is autonomous, the first part of
Theorem~\ref{teo:Pre:2.12} reduces to Theorem~2.12 in
\cite{rui-leite02-2}.
\end{remark}

Lemma~\ref{lemaRCDT} and Theorem~\ref{teo:Pre:2.12} give the main
motivation for our definition of generalized time-dependent spline
in $\mathbb{R}^n$. Let $L$ be the linear matrix differential
operator of order $p$
\begin{equation}
\label{eq:OopMatL} L = D^p \cdot - \,A_{p-1}(t)D^{p-1}\cdot -
\cdots - A_1(t)D\cdot - A_0(t)\cdot \, ,
\end{equation}
where each $A_j(t)$, $j = 0,1,\ldots,p-1$, is a real square
$n\times n$ $C^{\,p}$-smooth matrix function in $[a,b]$. The
operator $L$ is acting on the space $C^{\,m}[a,b]$ of real vector
functions defined in $[a,b]$. The adjoint of $L$, denoted by
$L^*$, is defined as
\[\begin{split}
L^* &= (-1)^p D^p\cdot +
\,(-1)^{p}D^{p-1}(A_{p-1}^{\prime}(t)\cdot) +
\,(-1)^{p-1}D^{p-2}(A_{p-2}^{\prime}(t)\cdot) + \\
&\quad +\cdots+ D(A_1^{\prime}(t)\cdot)-A_0^{\prime}(t)\cdot \, .
\end{split}\]
$L^*$ is also acting on $C^{\,m}[a,b]$ and the scalar product for
which it is computed is given by
\[\langle x_1 , x_2 \rangle = \int_{a}^{b} x_1(t)'x_2(t) \,dt\,.\]
Consider
\begin{equation}
\label{eq:Odelta} \Delta : a = t_0 < t_1 < \ldots < t_m = b
\end{equation}
to be a partition of $[a,b]$, and let $\Omega$ represent the set
of all $\mathbb{R}^n$-valued functions defined in $[a,b]$ which
are of class $C^{\,2p-2}$ in $[a,b]$ and of class $C^{\,2p}$ in
each interval $[t_i,t_{i+1}]$, $i=0,1,\ldots,m-1$.

\begin{definition}[Generalized time-dependent spline in $\mathbb{R}^n$]
\label{def:VGS:TD} A function $s : [a,b] \rightarrow \mathbb{R}^n$
is an interpolating generalized spline of $f\in\Omega$, associated
to $\Delta$ \eqref{eq:Odelta} and $L$ \eqref{eq:OopMatL}, if
$s\in\Omega$, $s$ is a solution of the matrix differential
equation $L^*Lx = 0$ in each interval $[t_i,t_{i+1}]$, $s(t)=f(t)$
on $\Delta$ (interpolation conditions), and
$s^{(k)}(t_0)=f^{(k)}(t_0)$, $s^{(k)}(t_m)=f^{(k)}(t_m)$, for
$k=1,2,\ldots,p-1$ (boundary conditions).
\end{definition}

\begin{remark}
Definition~\ref{def:VGS:TD} includes, as particular cases, the
scalar Definition~\ref{def:sg-tipoI} and the definition introduced
in \cite{rui-leite02-2}.
\end{remark}

\begin{remark}
As done in the scalar case, the interpolating function
$f\in\Omega$ can be omitted in Definition~\ref{def:VGS:TD}.
\end{remark}

\begin{remark}
The function $x(t)$, $t \in [a,b]$, given in each interval
$[t_{i},t_{i+1}]$, $i=0,1,\ldots,m-1$, by \eqref{eq:ExpExpGRSL},
is a generalized time-dependent spline in $\mathbb{R}^n$ in the
sense of Definition~\ref{def:VGS:TD}.
\end{remark}

\begin{remark}
For $L=D^p$ the solutions of $L^*Lx = 0$ give polynomial splines
in $\mathbb{R}^n$ with all the components being scalar polynomial
splines of degree $2p-1$. This is, as mentioned at the end of
\S\ref{subsec:sgs}, the immediate generalization of scalar
polynomial splines to vector-valued splines, and the one found in
the literature.
\end{remark}

We have seen that generalized splines associated to an operator
$L$ of order $p = 1$ are related to the optimal control problem
$(P)$. For $p > 1$, there corresponds an optimal control problem
with higher-order dynamic $x^{(p)} = \sum_{j=0}^{p-1}
A_j(t)x^{(j)} + B(t)u$. This higher-order optimal control problem
can be easily written in form $(P)$. For that we introduce new
state variables, reducing the control system of order $p$ to a
first-order control system. This is the same to say that when $L$
is an operator of order $p > 1$, the homogeneous differential
equation $L^*Lx = 0$ of order $2p$ can be reduced to a first order
differential equation, just by increasing the dimension of the
matrices $A_j(t)$, $j = 0,1,\ldots,p-1$.

Under our hypotheses, it is possible to write the optimal control
problem $(P)$ as a problem of the calculus of variations with
higher-order derivatives. This is done by showing that an
arbitrary admissible pair $\left(x(\cdot),u(\cdot)\right)$ of
$(P)$ can be always expressed in terms of higher order derivatives
of a single vector valued function (see \cite{MR91b:49033}). From
Theorem~\ref{teo:Pre:2.12} we obtain:

\medskip

\begin{theorem}
Given the operator $L$ \eqref{eq:OopMatL} and the partition $\Delta$
\eqref{eq:Odelta}, there exists a unique generalized spline in
$\mathbb{R}^n$ for each set of boundary and interpolation
conditions. This generalized spline is the unique solution of the
following higher-order problem of the calculus of variations:
\[\int_{a}^{b} \langle Lg,Lg \rangle \,dt \quad \rightarrow \quad \min \, ,\]
among all the functions $g\in\Omega$ that satisfy the same boundary
and interpolation conditions.
\end{theorem}


\section{Examples.}
We give two examples for which the state and control spaces are
$\mathbb{R}^2$. We denote the components of the state vector $x$
by $x_1$ and $x_2$; the components of the control vector $u$ by
$u_1$ and $u_2$. The first example is
\begin{equation*}
\min_{u \,=\, (u_1,u_2)^{\prime}}\quad
\int_{0}^{2}{u_1}^2+{u_2}^2\,dt
\end{equation*}
subject to the control system
\begin{equation}
\label{eq:ex1:NACS} \left\{\begin{array}{l}
\dot{x}_1 = t^2x_2 + u_2 \, ,\\
\dot{x}_2 = -t^2x_1 + u_1 \, ,
\end{array}\right.
\end{equation}
and the interpolating conditions
\begin{equation*}
x(t_0 = 0) = (0,0)^{\prime}, \quad x(t_1 = 1) = (1,0.5)^{\prime}
\, , \quad x(t_2 = 2) = (-0.25,1)^{\prime} \, .
\end{equation*}
As far as the control system \eqref{eq:ex1:NACS} is
non-autonomous, this example is not covered by the results in
\cite{rui-leite02-2}. We have the time interval $[0,2]$ and its
partition $\Delta : a = 0 < 1 < 2 = b$. The associated state
transition matrix is given by
\[\Phi(t,t_i) = \begin{pmatrix}
\cos{\left(\frac{t^3-{t_i}^3}{3}\right)} & \sin{\left(\frac{t^3-{t_i}^3}{3}\right)} \\
-\sin{\left(\frac{t^3-{t_i}^3}{3}\right)} &
\cos{\left(\frac{t^3-{t_i}^3}{3}\right)}
\end{pmatrix}.\]
The linear dynamic is completely state controllable. Such a
conclusion follows immediately from the fact that the symmetric
matrix
\[W=\int_{\tau_0}^{\tau_1}\Phi(\tau_0,s)B(s)B(s)^{\prime}\Phi(\tau_0,s)^{\prime}\,ds\]
is positive definite for some $\tau_1 > \tau_0$ with
$\tau_0,\tau_1\in[0,2]$. This is a classical test for complete
controllability which is due to Kalman \cite{kalman60}. Since $B$
and $\Phi$ are orthogonal matrices, the matrix $W$ is simply
\[\begin{pmatrix}
\tau_1-\tau_0 & 0 \\ 0 & \tau_1-\tau_0
\end{pmatrix}.\]
The optimal control is, in each interval $[t_i,t_{i+1}]$, solution
of the equation
\[L^*Bu = 0 \Leftrightarrow \left\{\begin{array}{l}
\dot{u}_1 + t^2u_2 = 0 \\
\dot{u}_2 - t^2u_1 = 0.
\end{array}\right.\]
We get
\begin{equation}u(t)_{|t\in[t_i,t_{i+1}]} = \left(\begin{array}{c}
-\sin{\left( \frac{t^3-{t_i}^3}{3} \right)}c_{1i}
+\cos{\left( \frac{t^3-{t_i}^3}{3} \right)}c_{2i} \\
\cos{\left(\frac{t^3-{t_i}^3}{3}\right)}c_{1i}
+\sin{\left(\frac{t^3-{t_i}^3}{3}\right)}c_{2i}\end{array}\right)
\label{eq:optimalcontrol2}\end{equation} where $c_{1i}$ and
$c_{2i}$ are real constants to be found. The corresponding
generalized spline in $\mathbb{R}^2$, solution of equation
\[L^*Lx = 0 \Leftrightarrow
\left\{\begin{array}{l}
\ddot{x}_1 - 2t^2\dot{x}_2 - t^4x_1 - 2tx_2 = 0 \\
\ddot{x}_2 + 2t^2\dot{x}_1 + 2tx_1 - t^4x_2 = 0
\end{array}\right.\]
in each interval $[t_i,t_{i+1}]$, is given by
$x(t)_{|t\in[t_i,t_{i+1}]} = (x_1(t), x_2(t))^{\prime}$ with
\[\textstyle x_1(t) = \cos{\left( \frac{t^3-{t_i}^3}{3}
\right)}(x_{1}(t_i)+(t-t_i)\,c_{1i} ) + \sin{\left(
\frac{t^3-{t_i}^3}{3}\right)}(x_{2}(t_i) + (t - t_i)\,c_{2i})\]
and
\[\textstyle x_2(t) = -\sin{\left(\frac{t^3-{t_i}^3}{3}\right)}(x_{1}(t_i)+(t-t_i)\,c_{1i})+
\cos{\left(\frac{t^3-{t_i}^3}{3}\right)}(x_{2}(t_i)+(t-t_i)\,c_{2i})\]
where $c_{1i}$ and $c_{2i}$ are the same constants that appear in
formula (\ref{eq:optimalcontrol2}). As expected, the resulting
spline is a continuous vector function and the optimal control
function is discontinuous at $t=t_1$.

\begin{center} \resizebox*{9cm}{5cm}{\includegraphics{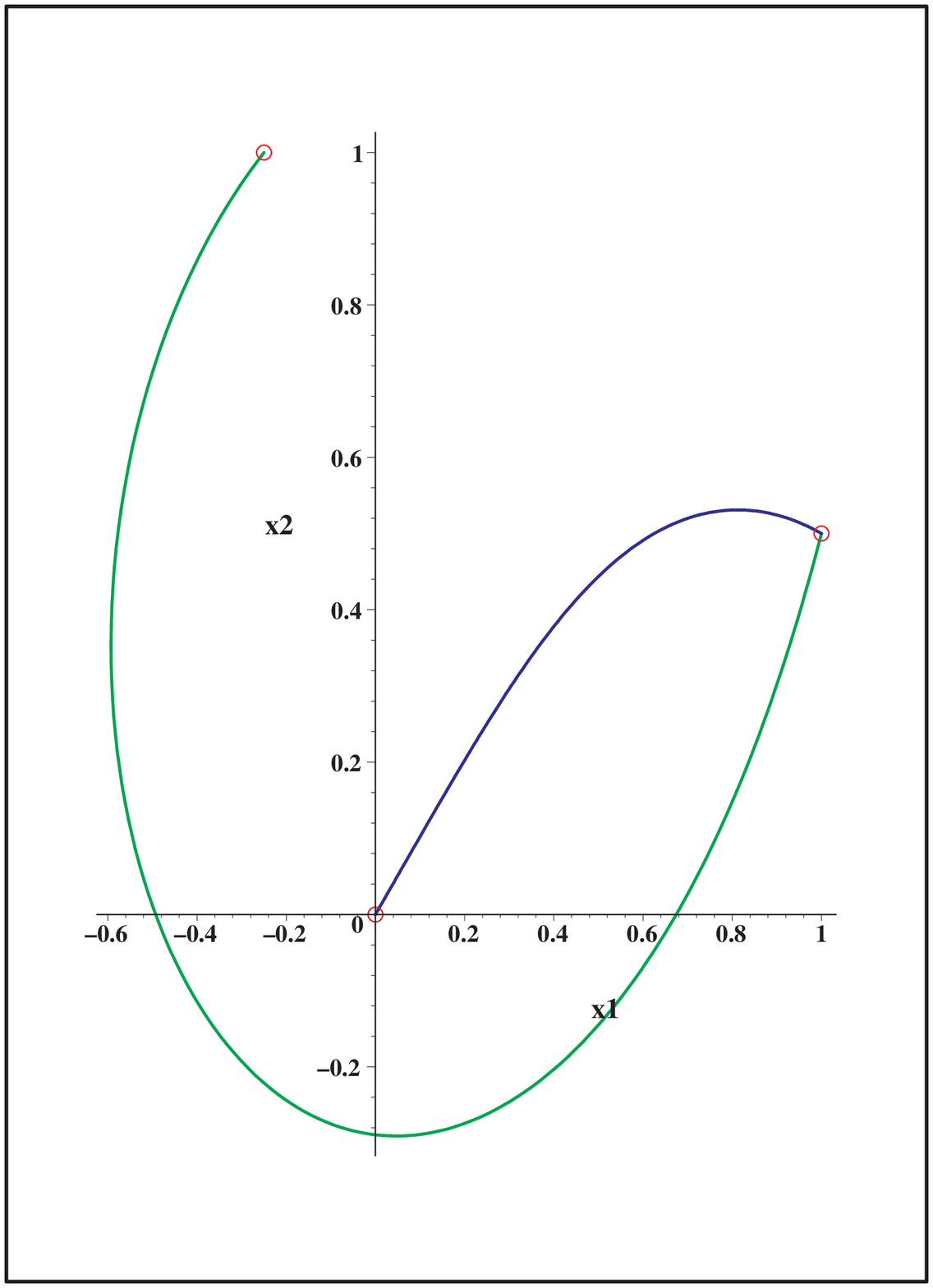}} \\
\mbox{} \\ First example -- generalized spline in $\mathbb{R}^2$
\end{center}

\begin{center} \resizebox*{9cm}{5cm}{\includegraphics{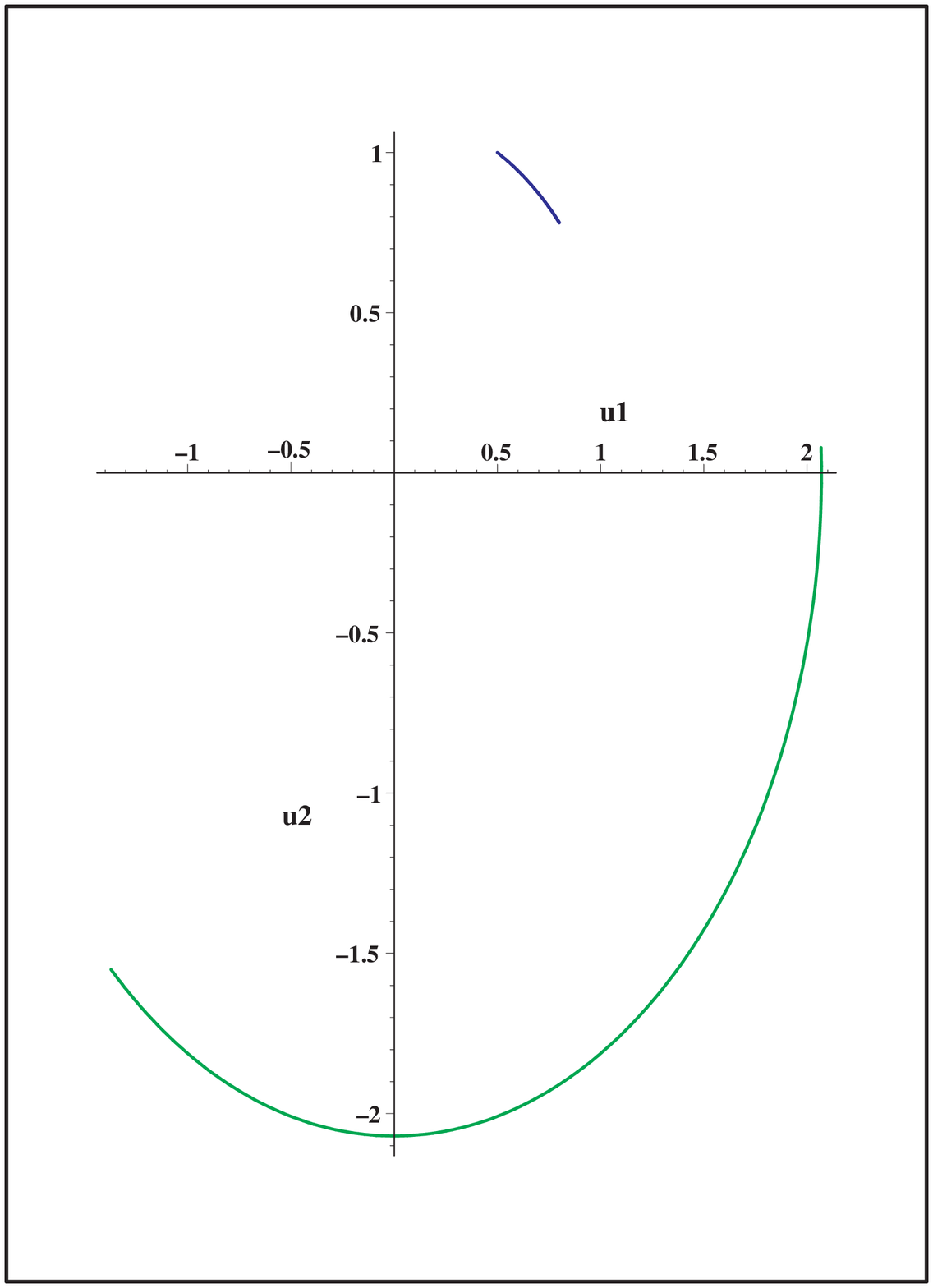}} \\
\mbox{} \\ First example -- optimal control
\end{center}

We now apply our results to the autonomous situation treated in
\cite{rui-leite02-2}. Consider the optimal control problem
\[\begin{split}
& \min_{u \,=\,(u_1,u_2)^{\prime}}\quad\int_{0}^{4}{u_1}^2+2u_1u_2
+ 2{u_2}^2\,dt \\
\text{subject to} \\
& \left\{\begin{array}{l}
\dot{x}_1 = -x_2 + u_2 \\
\dot{x}_2 = 2x_1 + u_1 + u_2,
\end{array}\right. \\
& x(t_0 = 0) = (0,0)^{\prime}, \quad x(t_1 = 1) = (1,0.5)^{\prime}, \\
& x(t_2 = 2) = (-0.25,1)^{\prime} \quad\text{and}\quad
x(t_3=4)=(1,-1)^{\prime}.
\end{split}\]
The optimal state trajectory is the generalized spline which, in
each interval $[t_i,t_{i+1}]$, is solution of equation
\[L^*Lx = 0 \Leftrightarrow \ddot{x} + (A^{\prime} - A)\,\dot{x}
- (A^{\prime}A)\,x=0 \Leftrightarrow \left\{\begin{array}{l}
\ddot{x}_1 + 3\dot{x}_2 - 4x_1 = 0 \\
\ddot{x}_2 - 3\dot{x}_1 - x_2 = 0.
\end{array}\right.\]
We get, in each interval $[t_i,t_{i+1}]$,
\[\textstyle x_1(t) = \sin{(\sqrt{2}\,t)}(\frac{3t}{4}c_{1i} + c_{4i})
+ \cos{(\sqrt{2}\,t)}( \frac{3\sqrt{2}}{8}c_{1i} +
\frac{3t}{4}c_{2i} + c_{3i})\]
and
\[\textstyle x_2(t) = \sin{(\sqrt{2}\,t)}( c_{1i} + \frac{3\sqrt{2}\,t}{4}c_{2i} +
\sqrt{2}\,c_{3i})+\cos{(\sqrt{2}\,t)}(-\frac{3\sqrt{2}\,t}{4}c_{1i}
+\frac{1}{4}c_{2i}+\sqrt{2}\,c_{4i})\]
where $c_{1i}$, $c_{2i}$,
$c_{3i}$ and $c_{4i}$ are real constants to be found in each
interval.

\begin{center} \resizebox*{9cm}{5cm}{\includegraphics{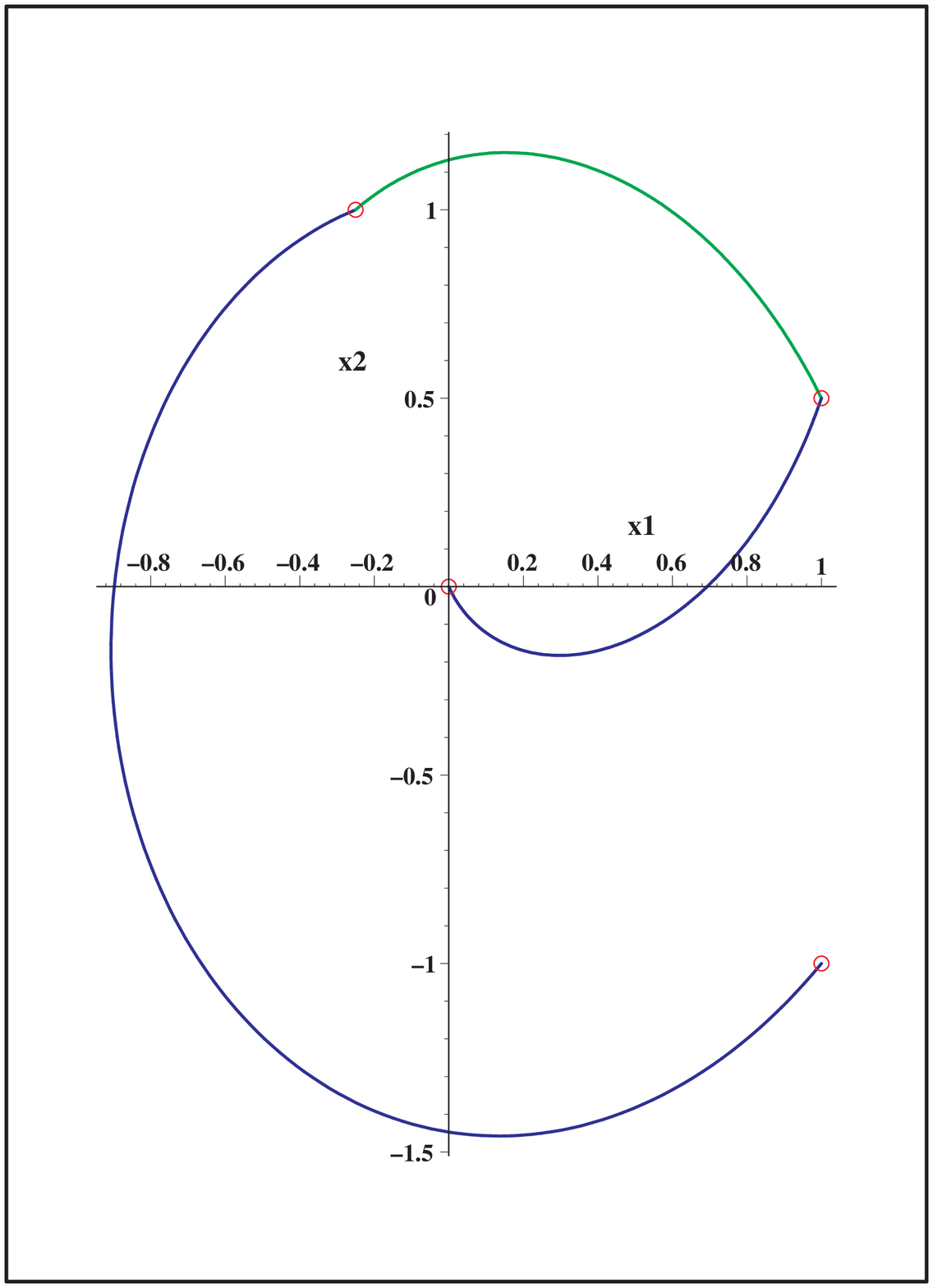}} \\
\mbox{} \\ Second example -- generalized spline in $\mathbb{R}^2$
\end{center}


\subsection*{Acknowledgements.}
R.~Rodrigues was supported in part by ISR-Coimbra and project
POSI/SRI/41618/2001. D.~Torres was supported by the R\&D unit CEOC
of the University of Aveiro, through the program POCTI of the
Portuguese Foundation for Science and Technology (FCT), cofinanced
by the European Community fund FEDER. The authors are grateful to
A. Sarychev who pointed out an inconsistency in an earlier version
of this paper.


\bigskip

\begin{flushleft}

Rui C. RODRIGUES\\
Department of Physics and Mathematics\\
Instituto Superior de Engenharia de Coimbra\\
Rua Pedro Nunes\\
3030-199 Coimbra, PORTUGAL\\
e-mail: \texttt{ruicr@isec.pt}\\[2ex]

Delfim F. M. TORRES\\
Department of Mathematics\\
University of Aveiro\\
3810-193 Aveiro, PORTUGAL\\
e-mail: \texttt{delfim@mat.ua.pt}\\[2ex]

\end{flushleft}


\begin{thebibliography}{10}

\bibitem{ahlberg67}\textsc{Ahlberg~J.~H., Nilson~E.~N., Walsh~J.~L.},
{\it The Theory of Splines and Their Applications}, Academic
Press, New York, 1967.

\bibitem{MR86k:49002} \textsc{Ball~J.~M. and Mizel~V.~J.},
\textit{One-dimensional variational problems whose minimizers do
not satisfy the {E}uler {L}agrange equation}, Arch. Rational Mech.
Anal. {\bf 90} 4 (1985), 325--388.

\bibitem{brockett70}\textsc{Brockett~R.~W.},
{\it Finite Dimensional Linear Systems}, John Wiley \& Sons, 1970.

\bibitem{MR85c:49001}\textsc{Cesari~L.},
{\it Optimization---theory and applications}, Springer-Verlag, New
York, 1983.

\bibitem{MR85m:49051} \textsc{Clarke~F.~H. and Vinter~R.~B.},
\textit{On the conditions under which the {E}uler equation or the
maximum principle hold}, Appl. Math. Optim. {\bf 12} 1 (1984),
73--79.

\bibitem{MR91b:49033} \textsc{Clarke~F.~H. and Vinter~R.~B.},
\textit{Regularity properties of optimal controls}, SIAM J.
Control Optim. {\bf 28} 4 (1990), 980--997.

\bibitem{MR22:13373} \textsc{Filippov~A.~F.},
\textit{On some questions in the theory of optimal regulation:
existence of a solution of the problem of optimal regulation in
the class of bounded measurable functions}, Vestnik Moskov. Univ.
Ser. Mat. Meh. Astr. Fiz. Him. 2 (1959), 25--32.

\bibitem{kalman60} \textsc{Kalman~R.~E.},
\textit{Contributions to the theory of optimal control}, Bol. Soc.
Mat. Mex. (1960), 102--119.

\bibitem{clyde95} \textsc{Martin~C., Enqvist~P., Tomlinson~J., Zhang~Z.},
\textit{Linear control theory, splines and interpolation},
Computation and Control {\bf iv} (1995), 269--287.

\bibitem{pontryagin62}\textsc{Pontryagin~L.~S., Boltyanskii~V.~G., Gamkrelidze~R.~V.,
Mishchenko~E.~F.}, {\it L. S. Pontryagin Selected Works - volume 4
- The Mathematical Theory of Optimal Processes}, Gordon and Breach
Science Publishers, 1986.

\bibitem{rui-leite02-2} \textsc{Rodrigues~R.~C. and Silva Leite~F.},
\textit{A multi-input/multi-output system representation of
generalized splines in $\mathbb{R}^n$}, Preprint, Department of
Mathematics, University of Coimbra 02-13 (2002).

\bibitem{rui-leite-simoes99} \textsc{Rodrigues~R.~C., Silva Leite~F., Sim\~{o}es~C.},
\textit{Generalized splines and optimal control}, Proceedings of
the European Control Conference, ECC'99, Karlsruhe, Germany,
(1999), CD-ROM paper \texttt{F0259.pdf}.

\bibitem{delfim03} \textsc{Torres~D.~F.~M.},
\textit{Lipschitzian regularity of the minimizing trajectories for
nonlinear optimal control problems}, Math. Control Signals Systems
{\bf 16} 2-3 (2003), 158--174.

\end{thebibliography}
\end{document}